\documentclass[11pt,a4paper]{amsart}
\usepackage{amsmath,amssymb}
\usepackage{amssymb,latexsym,bbm}
\usepackage{array}
\usepackage{mathrsfs}

\usepackage{amsthm}

{\theoremstyle{plain}
\newtheorem{definition}{Definition}[section]

\newtheorem{theorem}[definition]{Theorem}

\newtheorem{corr}[definition]{Corollary}

} {\theoremstyle{definition}
\newtheorem{ex}[definition]{Example}
\newtheorem{rem}[definition]{Remark}
}

 \DeclareMathOperator{\F}{{\mathcal F}}

\renewcommand{\mathbb}{\mathbbm}           
 \DeclareMathOperator{\C}{\mathbb{C}}
\DeclareMathOperator{\R}{\mathbb{R}} 
\DeclareMathOperator{\Rp}{\R_+} 

\renewcommand{\O}{{\mathcal O}}                      

\renewcommand{\epsilon}{\varepsilon}
\renewcommand{\phi}{\varphi}
 
\renewcommand{\theta}{\vartheta}
\renewcommand{\le}{\leqslant}
\renewcommand{\ge}{\geqslant}

\newcommand{\abs}[1]{\left\lvert #1 \right\rvert}   
\newcommand{\norm}[1]{\left\lVert #1 \right\rVert}  

\DeclareMathOperator{\1}{\mathbbm 1}                 

\DeclareMathOperator{\init}{\phi}

\begin{document}

\title[Geometric Brownian Motion with delay]
{Geometric Brownian Motion with delay:\\ mean square characterisation}
\author{John A. D. Appleby}
\address{School of Mathematical Sciences, Dublin City University, Dublin 9,
Ireland} \email{john.appleby@dcu.ie} \urladdr{http://webpages.dcu.ie/\textasciitilde
applebyj}
\author{Xuerong Mao}
\address{Department of Statistical and Modelling Science, Strathclyde
University, Glasgow, United Kingdom. } \email{xuerong@stams.strath.ac.uk}
\urladdr{http://www.stams.strath.ac.uk/\textasciitilde xuerong}
\author{Markus Riedle}
\address{Humboldt-University of Berlin,
  Department of Mathematics, Unter den Linden 6,
  10099 Berlin, Germany}
\email{riedle@mathematik.hu-berlin.de}
\urladdr{http://www.mathematik.hu-berlin.de/\textasciitilde riedle}

\thanks{John Appleby was partially funded by an Albert College Fellowship, awarded by Dublin City University's
Research Advisory Panel.}



\begin{abstract}
A geometric Brownian motion with delay is the solution of a
stochastic differential equation where the drift and diffusion
coefficient depend linearly on the past of the solution, i.e. a
linear stochastic functional differential equation. In this work the
asymptotic behavior in mean square of a geometric Brownian motion
with delay is completely characterized by a sufficient and necessary
condition in terms of the drift and diffusion coefficients.
\end{abstract}

\maketitle

\section{Introduction}
Geometric Brownian motion is one of the stochastic processes most
often used in applications, not least in financial mathematics for
modelling the dynamics of security prices. More recently, however,
modelling the price process by a geometric Brownian motion has
been criticised because the past of the volatility is not taken
into account.

The {\em geometric Brownian motion} is the strong solution of the
stochastic differential equation
\begin{align*}
 dX(t)=bX(t)\,dt + \sigma X(t)\,dW(t)\qquad\text{for } t\ge 0,
\end{align*}
where $b$ and $\sigma$ are some real constants. If we wish that
the dynamics of the process $X$ at time $t$ are to depend on its
past, a natural generalisation involves replacing the constants
$b$ and $\sigma$ by some linear functionals on an appropriate
function space, say the space of continuous functions on a bounded
interval. Then we are led by the Riesz' representation theorem to
the following stochastic differential equation with delay:
\begin{align}\label{eq.introstoch}
 dX(t)=\left(\int_{[-\alpha,0]}X(t+u)\mu(du)\right)dt
 + \left(\int_{[-\alpha,0]}X(t+u)\nu(du)\right)dW(t)
\end{align}
for all $t\ge 0$ and some measures $\mu,\,\nu$. We call the
solution $X$ of this stochastic differential equation {\em
geometric Brownian motion with delay} and its asymptotic behavior
in mean square will be characterised in this article. In contrast
to the geometric Brownian motion without delay no explicit
representation of $X$ is known from which the asymptotic behaviour
can be inferred directly.

Equation \eqref{eq.introstoch} is a stochastic functional
differential equation with diffusion coefficient depending on the
past. Such equations may exhibit most irregular asymptotic
behaviour, see for example Mohammed and Scheutzow~\cite{MoSch97}
for the noisy feedback equation. Nevertheless, for a much more
general equation than \eqref{eq.introstoch} a wide variety of
sufficient conditions have been established guaranteeing stability
in some sense. An exhaustive list of researchers and papers are
not quoted here, but a good selection of such results are collated
in Mao~\cite{Mao97} and Kolmanovskii and Myskhis~\cite{KolMys:99}.
Despite this activity over the last twenty--five years to the best
of our knowledge no sufficient and necessary conditions are known
for stability, even for the linear equation \eqref{eq.introstoch}.
By contrast, in this work we are able to find necessary and
sufficient conditions which characterise completely the asymptotic
behavior of the solution in mean square. An interesting by-product
of this stability characterisation is the observation that a
deterministic solution may transpire to be a solution of the
stochastic equation. This feature cannot arise in linear non-delay
stochastic equations.

The proof of our stability characterisation relies on the fact
that a non--negative functional of the process has expected value
which satisfies a deterministic linear renewal equation. The
asymptotic behaviour of this functional is characterised by the
renewal theorem; once this characterisation has been obtained, it
is straightforward to characterise the asymptotic behaviour of the
mean square.

\section{Preliminaries}

We first turn our attention to the deterministic delay equation
underlying the stochastic differential equation
\eqref{eq.introstoch}. For a fixed constant $\alpha\ge 0$ we
consider the deterministic linear delay differential equation
\begin{align}
  \begin{split}
     \dot{x}(t) &= \int_{[-\alpha,0]}x(t+u)\,\mu(du)   \quad\text{for }t\ge 0,\\
     x(t) &= \init(t) \quad \text{for }t\in [-\alpha,0],
  \end{split}\label{eq.det}
\end{align}
for a measure $\mu\in M= M[-\alpha,0]$, the space of signed Borel
measure on $[-\alpha,0]$ with the total variation norm
$\norm{\cdot}_{TV}$. The initial function $\init$ is assumed to be
in the space $C[-\alpha,0]:=\{\psi:[-\alpha,0]\to\R:
\text{continuous}\}$. A function $x:[-\alpha,\infty)\to\R$ is
called a {\em solution} of \eqref{eq.det} if $x$ is continuous on
$[-\alpha,\infty)$, its restriction to $[0,\infty)$ is
continuously differentiable, and $x$ satisfies the first and
second identity of \eqref{eq.det} for all $t\ge 0$ and $t\in
[-\alpha,0]$, respectively. It is well known that for every
$\init\in C[-\alpha,0]$ the problem \eqref{eq.det} admits a unique
solution $x = x(\cdot,\init)$.

The {\em fundamental solution} or {\em resolvent} of \eqref{eq.det} is the unique locally
absolutely continuous function $r: [0,\infty)\to\R$ which satisfies
$$ r(t) = 1+ \int_0^t \int_{[\max\{-\alpha,-s\},0]} r(s+u)\,d\mu(u)\,ds \quad\text{for} \ t\ge 0.$$
It plays a role which is analogous to the fundamental system in
linear ordinary differential equations and the Green function in
partial differential equations. Formally, it is the solution of
\eqref{eq.det} corresponding to the initial function $\init =
\1_{\{0\}}$. For later convenience we set $r(u)=0$ for $u\in
[-\alpha,0)$.

The solution $x(\cdot,\phi)$ of \eqref{eq.det} for
an arbitrary initial segment $\init$ exists, is unique, and can be represented as
\begin{equation}\label{eq.xdet}
x(t,\init)=\phi(0)r(t)+\int_{[-\alpha,0]} \int_{s}^0 r(t+s-u)\init(u)\,du\,\mu(ds)
\quad\text{for }t\ge 0,
\end{equation}
cf. Diekmann et al \cite[Chapter I]{Dieetal95}. The fundamental
solution converges for $t\to\infty$ to zero if and only if
\begin{align}\label{eq.v_0}
 v_0(\mu)
 :=\sup\left\{\Re(\lambda):\, \lambda\in\C,\;\lambda - \int_{[-\alpha,0]}e^{\lambda s} \,
\mu(ds)=0\right\}
 <0,
\end{align}
where $\Re(z)$ denotes the real part of a complex number $z$. In
this case the decay is exponentially fast (see Diekmann et al
\cite[Thm. 5.4]{Dieetal95}) and the zero solution of
\eqref{eq.det} is uniformly asymptotically stable. In this
situation we have for every solution $x(\cdot,\phi)$ of
\eqref{eq.det}:
\begin{align*}
 \abs{x(t,\phi)}\le c e^{v_0(\mu)t}\qquad \text{for all }t\ge 0,
\end{align*}
and for a constant $c>0$ depending only on the initial function
$\phi$.

Let us introduce equivalent notation for \eqref{eq.det}. For a
function $x:[-\alpha,\infty)\to \R$ we define the {\em segment of
$x$} at time $t\ge 0$ by the function
\begin{align*}
x_t:[-\alpha,0]\to \R, \qquad x_t(u):=x(t+u).
\end{align*}
If we equip the space $C[-\alpha,0]$ of continuous functions with
the supremum norm Riesz' representation theorem guarantees that
every continuous functional $F:C[-\alpha,0]\to\R$ is of the form
\begin{align*}
 F(\psi)=\int_{[-\alpha,0]}\psi(u)\,\mu(du)
\end{align*}
for a measure $\mu\in M$. Hence, we will write \eqref{eq.det} in
the form
\begin{align*}
 \dot{x}(t)=F(x_t)\quad \text{for }t\ge 0,\qquad
     x_0 &= \init
\end{align*}
and assume $F$ to be a continuous and linear functional on
$C[-\alpha,0]$.

Let us fix a complete probability space $(\Omega,\F,P)$ with a filtration $(\F_t)_{t\ge
0}$ satisfying the usual conditions and let $(W(t):\,t\ge 0)$ be a Wiener process on this
space. We study the following stochastic differential equation with time delay:
\begin{align}
\begin{split}
 dX(t)&= F(X_t)\,dt+ G(X_t)\,dW(t)\quad\text{for }t\ge 0,\\
  X(u)&= \init(u)\quad\text{for } u\in [-\alpha,0],
 \end{split}\label{eq.stoch}
\end{align}
where $F$ and $G$ are continuous and linear functionals on
$C[-\alpha,0]$ for a constant $\alpha\ge 0$. As before, we can
write the functional $G$ in the form
\begin{align*}
 G(\psi)=\int_{[-\alpha,0]}\psi(u)\,\nu(du) \qquad\text{for all }
  \psi\in C[-\beta,0]
\end{align*}
and for a measure $\nu\in M$. We note that assuming the same
domain $[-\alpha,0]$ for the arguments of the functionals $F$ and
$G$ does not involve any restriction or loss of generality.

For every $\phi\in C[-\alpha,0]$ there exists a unique, adapted
strong solution $(X(t,\phi):\,t\ge -\alpha)$ with finite second
moments of \eqref{eq.stoch} (cf., e.g., Mao~\cite{Mao97}). The
dependence of the solutions on the initial condition $\phi$ is
neglected in our notation in what follows; that is, we will write
$x(t)=x(t,\phi)$ and $X(t)=X(t,\phi)$ for the solutions of
\eqref{eq.det} and \eqref{eq.stoch} respectively.

By Rei{\ss} et al \cite[Lemma 6.1]{ReRiGa06} the solution
$(X(t):\,t\ge -\alpha)$ of \eqref{eq.stoch} obeys a variation of
constants formula
\begin{align}\label{eq.varofconst}
 X(t)= \begin{cases}
 x(t)+ \int_0^t r(t-s) G(X_s)\,dW(s),& t\ge 0,\\
  \init(u) , & u\in [-\alpha,0],
  \end{cases}
\end{align}
where $r$ is the fundamental solution of \eqref{eq.det}. It is to
be noted that this equation does not supply an explicit form of
the solution.

\section{Stability}

The asymptotic behavior of the solution $X$ relies on the stochastic
convolution integral arising in the variation of constants formula
\eqref{eq.varofconst}. Let us define
\begin{align}\label{eq.defY}
 Y(t):= G(X_t)
 \qquad \text{for }t\ge 0,
\end{align}
such that the stochastic convolution integral is the convolution
of the stochastic process $Y=(Y(t):\,t\ge 0)$ and the fundamental
solution $r$. The following result shows that the functional
$\mathbb{E}[Y^2]$ satisfies a linear convolution integral
equation.
\begin{theorem}
Let $(X(t):\,t\ge -\alpha)$ be the solution  of \eqref{eq.stoch}. Then we have
for all $ t \ge 0$
\begin{align}\label{eq.EX}
 E\abs{X(t)}^2= \abs{x(t)}^2 + \int_0^t r^2(t-s)\,E\abs{Y(s)}^2 \,
 ds,
\end{align}
where $Y$, defined by \eqref{eq.defY}, obeys for all $t\ge 0$:
\begin{align}\label{eq.EY}
 E\abs{Y(t)}^2= G^2(x_t)+ \int_0^t G^2(r_{t-s}) E\abs{Y(s)}^2\, ds.
\end{align}
\end{theorem}

\begin{proof}
The variation of constants formula \eqref{eq.varofconst} and It{\^o}'s isometry imply the
first assertion.

Using again the variation of constants formula we obtain by
Fubinis's theorem for stochastic integrals for $t\in [0,\alpha]$:
\begin{align*}
 E\abs{Y(t)}^2& = E\abs{G(X_t)}^2\\
 &=E\abs{\int_{[-\alpha,-t)}X_t(u)\,\nu(du) +
 \int_{[-t,0]}X_t(u)\,\nu(du)}^2 \\
 &=E\abs{\int_{[-\alpha,-t)}\phi(t+u)\,\nu(du)+
  \int_{[-t,0]}x(t+u)\,\nu(du)\right.\\
&\qquad + \left.  \int_{[-t,0]}\left( \int_0^{t+u}
r(t+u-s)Y(s)\,dW(s)\,\right)
    \, \nu(du)}^2\\
& = E\abs{G(x_t) +  \int_0^t \left( \int_{[s-t,0]}
r(t+u-s)\,\nu(du)\right)
        Y(s)\, dW(s)}^2\\
 &= \abs{G(x_t)}^2 + \int_0^t G^2(r_{t-s}) E\abs{Y(s)}^2\, ds,
\end{align*}
where we used in the last line $r(u)=0$ for $u<0$. Setting
$\nu([a,b))=0$ for all $a\le b\le -\alpha$ enables us to enlarge the
integration domain for $G$ such that we can write also for $t\ge
\alpha$:
\begin{align*}
 E\abs{Y(t)}^2& =E
 \abs{G(x_t) +  \int_{[-t,0]}\left(
\int_0^{t+u} r(t+u-s)Y(s)\,dW(s)\,\right)
    \, \nu(du)}^2.
\end{align*}
We can proceed as above to verify that equation \eqref{eq.EY} is
also satisfied for $t\ge \alpha$.
\end{proof}

We next turn to stating and proving our first stability result. In
it, the hypothesis $r\in L^2(\Rp)$ is employed. We remark that
this assumption is necessary if $E\abs{X^2(t)}\to0$ as
$t\to\infty$. To see this, note by \eqref{eq.EX} that
$E\abs{X^2(t)}$ cannot tend to zero as $t\to\infty$ if $x(t)$ does
not tend to zero. But the latter cannot occur if $r$ is not in
$L^2(\Rp)$.

The function $s\mapsto G(r_s)$, which we denote by
$G(r_{\bullet})$, is square integrable if $r\in L^2(\Rp)$ and its
norm in $L^2(\Rp)$ is given by
\begin{align*}
 \norm{G(r_\bullet)}_{L^2(\Rp)}=\left(\int_0^\infty
 (G(r_s))^2\,ds\right)^{1/2}.
\end{align*}
This quantity allows to characterise the asymptotic behaviour of
the solution for \eqref{eq.stoch}:

\begin{theorem}\label{th.dich}
If the fundamental solution $r$ is in $L^2(\Rp)$ then the solution
$(X(t):\,t\ge -\alpha)$ of \eqref{eq.stoch} obeys the following
trichotomy:
\begin{enumerate}
\item[{\rm (a)}] if $\norm{G(r_{\bullet})}_{L^2(\Rp)} <1$, then
there exists  $\kappa>0$ such
 \begin{align*}
 \lim_{t\to\infty}e^{\kappa t} E\abs{X(t)}^2=0.
\end{align*}
\item[{\rm (b)}] if $\norm{G(r_{\bullet})}_{L^2(\Rp)} =1$, then
\begin{align*}
   \lim_{t\to\infty}E\abs{X(t)}^2 =
\frac{\left(\int_0^\infty G^2(x_s) \,ds\right)\left( \int_0^\infty
r^2(s)\,ds\right)}{\int_0^\infty s G^2(r_s)\,ds}<\infty.\\
\end{align*}
\item[{\rm (c)}] if $\norm{G(r_{\bullet})}_{L^2(\Rp)} >1$, then
there exists $\kappa>0$ such that
\begin{align*}
  \qquad \lim_{t\to\infty}e^{-\kappa t} E\abs{X(t)}^2 =
\frac{\left(\int_0^\infty e^{-\kappa s} G^2(x_s)\,ds
\right)\left(\int_0^\infty e^{-\kappa s}
  r^2(s)\,ds \right)}{\int_0^\infty se^{-\kappa s}G^2(r_s)\,ds}<\infty.\\
\end{align*}
\end{enumerate}
 \end{theorem}
\begin{proof}
Let us introduce the following functions and measures for $t\ge
0$:
\begin{align*}
 y(t):= E\abs{Y(t)}^2 \qquad f(t):= G^2(x_t)\\
 g(t):=G^2(r_t) \qquad \zeta(dt):= g(t)\, dt.
\end{align*}
Then we can rewrite \eqref{eq.EY} as the renewal equation
\begin{align}\label{eq.renewal}
 y(t)= f(t) + \int_0^t y(t-s)\,\zeta(ds) \quad\text{for all }t\ge 0
\end{align}
and the three cases (a) to (c) correspond to whether the renewal
equation \eqref{eq.renewal} is defective, proper or excessive.

To give the main idea of the proof we establish (a) first without
the convergence rate. In case (a) the renewal Theorem \cite[Thm
3.1.4]{Als91} implies
\begin{align*}
 \lim_{t\to\infty} y(t)= \frac{f(\infty)}{1-\norm{G(r_{\bullet})}_{L^2(\Rp)}}=0,
  \end{align*}
as $f(\infty):=\lim_{t\to \infty}f(t)=0$ due to $x(t)\to 0$ for
$t\to\infty$. Using the notation for $y$ introduced above,
equation \eqref{eq.EX} reads
\begin{align}\label{eq.EX2}
 E \abs{X(t)}^2 = x^2(t) + \int_0^t r^2(s)y(t-s)\,ds.
\end{align}
Consequently, as $x(t)\to 0$ for $t\to\infty$ and $r\in L^2(\Rp)$ we
arrive at
\begin{align*}
\lim_{t\to\infty} E\abs{X(t)}^2
 =\lim_{t\to\infty} \int_0^\infty r^2(s)\1_{[0,t]}(s) y(t-s)\,ds =0
\end{align*}
by dominated convergence.

We now turn to prove the exponential decay. Because
$\norm{G(r_{\bullet})}_{L^2(\Rp)} <1$ there exists $\theta>0$ such
that
\begin{align*}
 \int_0^\infty e^{\theta s}\, \zeta(ds)=1.
\end{align*}
Moreover, $r\in L^2(\Rp)$ implies that there exists $\gamma>0$
such that $r(t)=\O(\exp(-\gamma t))$ and $x(t)=\O(\exp(-\gamma
t))$ as $t\to\infty$. Consequently, we have $f(t)=\O(\exp(-2
\gamma t))$ and $g(t)=\O(\exp(-2\gamma t))$ as $t\to\infty$.
Therefore, we can infer by standard methods in renewal theory that
$y(t)=o(\exp(-\kappa t))$ for all $\kappa <(2\gamma\wedge\theta)$
which leads to
\begin{align*}
 \int_0^t r^2(s)y(t-s)\,ds =
  o\big(\exp(-\kappa t\big)).
\end{align*}
Consequently, equation \eqref{eq.EX2} yields the assertion for all
$\kappa < 2\gamma \wedge \theta$.

In case (b), the renewal  Theorem \cite[Thm 3.1.5]{Als91} implies
under appropriate conditions on $f$ that
\begin{align*}
 \lim_{t\to\infty} y(t)= (m(\zeta))^{-1} \int_0^\infty f(s)\,ds
  \end{align*}
with $m(\zeta):=\int_0^\infty s G^2(r_s)\,ds$. Note that
$m(\zeta)$ is finite as $r$ tends to zero exponentially fast and
it is non-zero because $\int_0^\infty G^2(r_s)\,ds=1$. Since the
measure $\zeta$ is absolutely continuous with respect to the
Lebesgue measure it is sufficient for the application of the
renewal Theorem that the function $f$ be in $(L^1\cap
L^\infty)(\Rp)$ and $f(t)\to 0$ for $t\to\infty$: both of these
conditions are satisfied here. By proceeding as in case (a) we
obtain
\begin{align*}
\lim_{t\to\infty} E\abs{X(t)}^2
  &= \frac{1}{m(\zeta)}\left(\int_0^\infty f(s) \,ds\right) \left(\int_0^\infty r^2(s)\,ds
  \right).
\end{align*}

In case (c), the renewal equation \eqref{eq.renewal} is excessive.
Then there exists a unique $\kappa>0$ such that
$\int_{\Rp}e^{-\kappa s}\,\zeta(ds)=1$. Furthermore,
$m_\kappa(\zeta):=\int_0^\infty se^{-\kappa s}\,\zeta(ds)$ is
non-zero (see \cite[Remark 3.1.8]{Als91}) and $m_\kappa(\zeta)$ is
finite as $r$ decays exponentially fast. As in case (b) it is
sufficient for the application of the renewal Theorem that the
function $f_\kappa$ with $f_\kappa(t):=e^{-\kappa t}f(t)$ for
$t\ge 0$ tends to zero and is in $L^1(\Rp)$. These conditions are
satisfied as $x(\cdot,\phi)$ decays exponentially fast and hence,
the renewal Theorem \cite[Cor. 3.1.9]{Als91} implies
\begin{align*}
   \lim_{t\to\infty}e^{-\kappa t } y(t)
   = (m_\kappa(\zeta))^{-1} \int_0^\infty e^{-\kappa s}f(s)\,ds.
\end{align*}
Finally, from \eqref{eq.EX} we have
\begin{align*}
 e^{-\kappa t} E\abs{X(t)}^2
  = e^{-\kappa t}\abs{x(t)}^2 + \int_0^t e^{-\kappa s} r^2(s)
     e^{-\kappa(t-s)}y(t-s)\,ds,
\end{align*}
and so, because of the exponential decay of $r$, we conclude
\begin{align*}
 \lim_{t\to\infty} e^{-\kappa t}E\abs{X(t)}^2
  = \frac{1}{m_\kappa (\zeta)}\left(\int_0^\infty e^{-\kappa s} f(s)\,ds
  \right)\left(\int_0^\infty e^{-\kappa s}
  r^2(s)\,ds \right).
\end{align*}
\end{proof}
Alsmeyer~\cite{Als91} contains a treatment of the renewal equation
which covers equations with measures. Similar results may be found
in Feller~\cite{Feller:1}.

As a corollary of Theorem~\ref{th.dich} we obtain an equivalence
between the asymptotic behavior of the mean square of the solution
$X$ and a condition on the fundamental solution $r$ and the
diffusion coefficient $G$. Naturally, this requires that the
solution $X$ does not reduce to the solution of the deterministic
equation \eqref{eq.det}, for in this case $X$ would not provide
any information on the diffusion coefficient. We argue below that
this situation may occur, and must be excluded in the next
corollary. This will also be illustrated presently in an example.
\begin{corr}\label{co.dich}
Let the fundamental solution $r$ be in $L^2(\Rp)$ and assume that
no version of the solution $X=X(\cdot,\phi)$ of \eqref{eq.stoch}
coincides with the deterministic solution $x=x(\cdot,\phi)$ of
\eqref{eq.det}. Then we have the following:
\begin{align*}
\lim_{t\to\infty} E\abs{X(t)}^2 =\begin{cases}
 0     &\Longleftrightarrow\; \; \norm{G(r_{\bullet})}_{L^2(\Rp)} <1, \\
 c>0   &\Longleftrightarrow\; \; \norm{G(r_{\bullet})}_{L^2(\Rp)} =1, \\
 \infty &\Longleftrightarrow\; \; \norm{G(r_{\bullet})}_{L^2(\Rp)} >1.
\end{cases}
\end{align*}
\end{corr}

\begin{proof}
We have to show that the constants in Theorem \ref{th.dich} in (b)
and (c) describing the limiting behavior of $X$ are non-zero, which
is equivalent to
\begin{align}\label{eq.intG=0}
 \int_0^\infty G^2(x_s(\cdot,\phi))\,ds \neq 0.
\end{align}
Because $s\mapsto x_s(\cdot,\phi)$ and $G$ are continuous
operators, equation \eqref{eq.intG=0} does not hold if and only if
$G(x_t(\cdot,\phi))=0$ for all $t\ge 0$. In this case
$x(\cdot,\phi)$ would be also a solution of the stochastic
equation \eqref{eq.stoch} which contradicts our assumption because
of the uniqueness of the solution of \eqref{eq.stoch}.
\end{proof}

\begin{rem} \label{rem.3.4}
If the solution $x(\cdot,\phi)$ of the deterministic equation
\eqref{eq.det} also solves the stochastic equation
\eqref{eq.stoch} it follows by taking expectation and It{\^o}'s
isometry that
\begin{align}\label{eq.remdich}
 G(x_t(\cdot,\phi))= 0\qquad\text{for all }t\ge 0.
\end{align}
Conversely, condition \eqref{eq.remdich} implies that a version of
the solution $X(\cdot,\phi)$ of \eqref{eq.stoch}  coincides with
$x(\cdot,\phi)$.

Thus, a sufficient condition for the trichotomy only relying on
the initial condition $\phi$ is $G(\phi)\neq 0$.

A more abstract point of view  shows that equation
\eqref{eq.remdich} holds true if a generalized eigenspace $N$ of
the deterministic equation \eqref{eq.det} is a subset of the
kernel ker$(G)$ of the diffusion coefficient $G$. Then, for every
$\phi\in N$ the segment $x_t(\cdot,\phi)$ is in $N$ and
consequently in ker$(G)$, so that \eqref{eq.remdich} is satisfied.
For details on eigenspaces and related results for equation
\eqref{eq.det} we refer the reader to Diekmann et al
\cite{Dieetal95} and Hale and Lunel \cite{HaleLun93}. A concrete
example is given below.

We emphasise that the forgoing situation in which a deterministic
solution may solve a non-trivial linear stochastic differential
equation is a very specific feature of stochastic functional
differential equations, and cannot occur in linear stochastic
ordinary differential equations.
\end{rem}

\begin{ex}
Let us consider the solution $(X(t):\,t\ge 0)$ of the  simple equation
\begin{align}\label{eq.ex}
 dX(t)=bX(t)\,dt + (cX(t)+dX(t-\alpha))\,dW(t)\qquad\text{for }t\ge
 0,
\end{align}
with $b<0$, $c,d\in\R$ and $\alpha >0$. For the condition in Theorem
\ref{th.dich} we calculate
\begin{align*}
 \norm{G(r_\bullet)}^2_{L^2(\Rp)}=\int_0^\infty G^2(r_s)\,ds
 = \frac{1}{-2b}(c^2+d^2+2cd e^{b\alpha}).
\end{align*}
By using results on deterministic linear difference equations we
obtain for a continuous function $y:[-\alpha,\infty)\to\R$:
\begin{align*}
 G(y_t)=0 \quad\text{ for all $t\ge 0$}
  & \;\;\;\Longleftrightarrow\;\;\;
 c y(t)+dy(t-\alpha)=0 \qquad\text{for all $t\ge 0$}\\
 & \;\;\;\Longleftrightarrow\;\;\;
 y(t)=y(0)e^{\gamma t} \quad\text{ for all $t\ge -\alpha$}
  \quad\text{and}\quad cd<0
\end{align*}
with $\gamma:=\tfrac{1}{\alpha}\ln \frac{-d}{c}$. In the case when
$b\neq \gamma$ and $cd<0$, we obtain that for every initial
condition $\phi$ that the solution $X=X(\cdot,\phi)$ obeys:
\begin{align}\label{eq.exstable}
 \lim_{t\to\infty} E\abs{X(t)}^2=0
 \;\;\Longleftrightarrow \;\; c^2+d^2+2cde^{b\alpha}<-2b.
\end{align}
In the case when $b<0$ and $cd\ge 0$, the equivalence
\eqref{eq.exstable} also holds.

In the non-delay case $d=0$ the solution $X$ is the geometric
Brownian motion and \eqref{eq.exstable} reproduces the well-known
and easily calculated fact that $E\abs{X(t)}^2\to 0$ if and only
if $c^2<-2b$. In the pure delay case $c=0$ and $d\neq 0$ we find
that $d^2<-2b$ is necessary and sufficient to guarantee
$E\abs{X(t)}^2\to 0$. However, although the condition on the noise
intensities is the same as for geometric Brownian motion, the rate
of decay to zero is different.

If $c\neq 0$ and $d\neq 0$ then the dependence of the stability
region on  the drift coefficient $b$ is described by $b < b_0<0$
where $b_0$ is the largest real root of
\begin{align*}
  c^2+d^2 + 2cde^{b_0\alpha} +2b_0=0.
\end{align*}
In particular, we observe that while the stability region for
\eqref{eq.ex} is symmetric in $c$ and $d$, it is not symmetric in
the sign of $cd$.

We finish by pointing out that equation \eqref{eq.ex} provides an
example of the situation already mentioned in Remark \ref{rem.3.4}
in which a solution of the deterministic equation \eqref{eq.det}
is also a solution of the stochastic equation \eqref{eq.stoch}. To
see this take for example $c=-e$, $d=1$ and $b=\gamma$ for some
$\alpha>0$. Then, for $\phi(u)=e^{\gamma u}$, $u\in [-\alpha,0]$,
the solution $x(t,\phi)=e^{\gamma t}$, $t\ge 0$, satisfies
$G(x_t(\cdot,\phi))=0$ for all $t\ge 0$. Thus, $x(\cdot,\phi)$ is
also a solution of the stochastic equation and $x(t,\phi)\to 0 $
for $t\to\infty$ even though $\norm{G(r_\bullet)}^2_{L^2(\Rp)}>1$.
\end{ex}

\end{document}